\documentclass[twoside,11pt]{article}
\usepackage{geometry}
\usepackage{graphicx}
\usepackage{indentfirst}                                % make the first line of all sections be indented
\usepackage[colorlinks]{hyperref}
\hypersetup{linkcolor=blue,filecolor=black,urlcolor=blue, citecolor=black}   % book-mark black
\usepackage[numbers,sort&compress]{natbib}         % references

\ifxetex
\usepackage{letltxmacro}
\LetLtxMacro\SavedIncludeGraphics\includegraphics
\def\includegraphics#1#{% #1 catches optional stuff (star/opt. arg.)
 \IncludeGraphicsAux{#1}%
}%
\newcommand*{\IncludeGraphicsAux}[2]{%
 \XeTeXLinkBox{%
  \SavedIncludeGraphics#1{#2}%
}}
\fi
\newcommand\orcidicon[1]{\href{https://orcid.org/#1}{\includegraphics[scale=0.02]{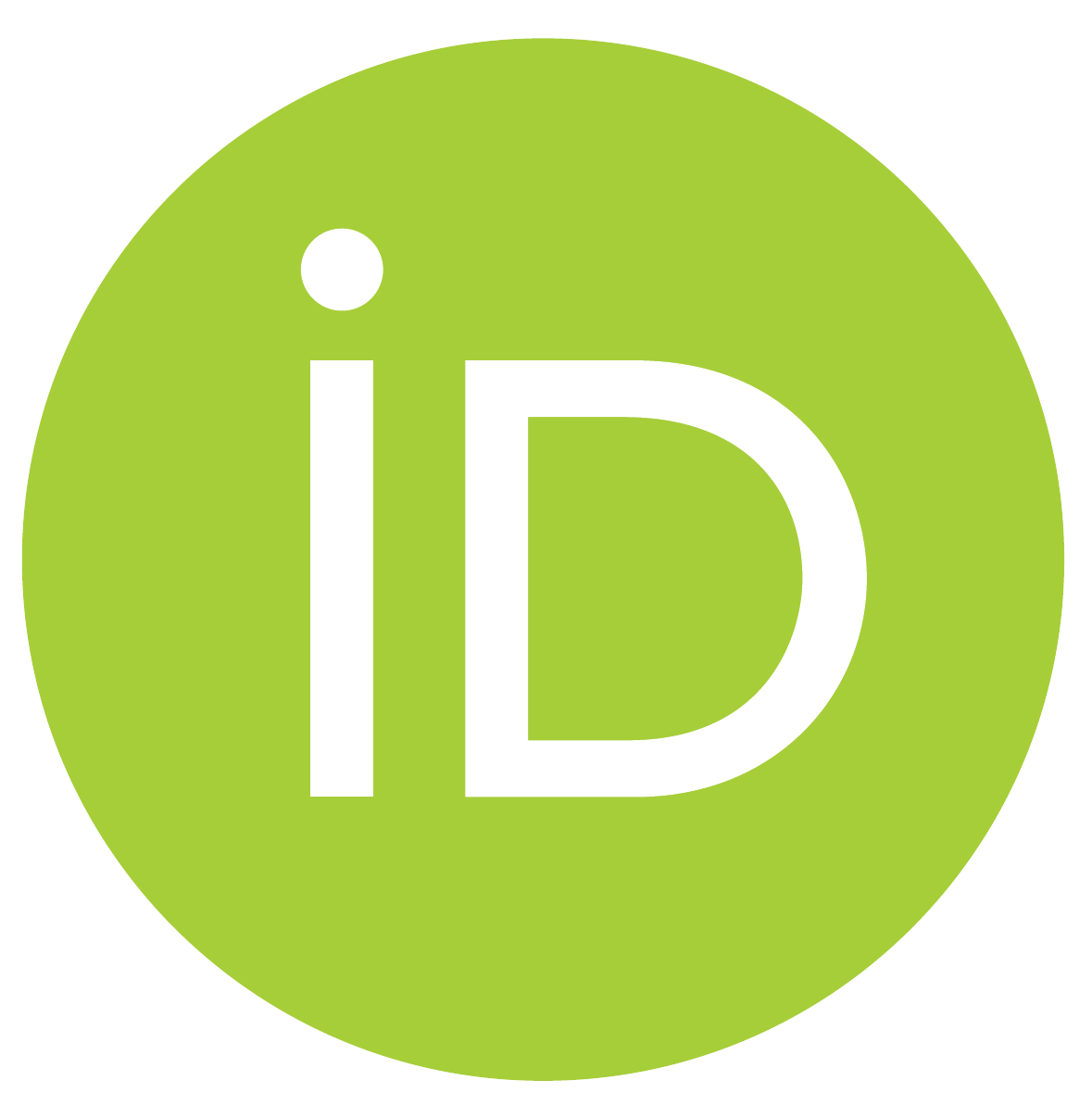}}}

\usepackage[titletoc,title]{appendix}
\usepackage{titlesec}
\titlespacing{\section}{0pt}{2.5ex}{1.5ex}
\titlespacing{\subsection}{0pt}{1.5ex}{1ex}
\titlespacing{\subsubsection}{0pt}{1.5ex}{1ex}
\titleformat{\section}{\large\bfseries\centering}{\thesection}{1em}{}
\titleformat{\subsection}[runin]{\bfseries}{\thesubsection.}{0.5em}{}[.\mbox{\ }]
\titleformat{\subsubsection}[runin]{\bfseries}{\thesubsubsection.}{0.4em}{}[.\mbox{\ }]

\usepackage{amsmath,amsfonts,amssymb,bm,mathrsfs,amsthm,amsxtra}
\numberwithin{equation}{section}
\newtheorem{lemma}{Lemma}[section]
\newtheorem{proposition}[lemma]{Proposition}
\newtheorem{theorem}{Theorem}[section]

\def\nt{|\hspace{-0.7pt}|\hspace{-0.7pt}|}

\usepackage{ulem}%%%% bar the text
\newcommand{\VERT}{\vert\kern-0.3ex\vert\kern-0.3ex\vert}
\def\p{\partial}

\def\div{{\rm div}\,}
\usepackage{accents}
\makeatletter
\def\widebar{\accentset{{\cc@style\underline{\mskip10mu}}}}
\makeatother

\newcommand\w[1]{\makebox[1em]{$#1$}}

\begin{document}
%%%%%----------
%%%%%---------- Full Title
%%%%%---------- Vacuum Free Boundary Problems in Ideal Compressible MHD
\title{\bf Stabilizing effect of surface tension \\ for the linearized MHD--Maxwell \\ free interface problem\let\thefootnote\relax\footnotetext{
This research was carried out at the Sobolev Institute of Mathematics, under a state contract (project no. FWNF-2022-0008).
}
}

%%%%%----------
%%%%%---------- Full Author
%%%%%----------
\author{
{\sc Yuri Trakhinin}\orcidicon{0000-0001-8827-2630}\thanks{Sobolev Institute of Mathematics, Koptyug av.~4, 630090 Novosibirsk, Russia. e-mail: trakhin@math.nsc.ru}
}

\date{ }

\maketitle

\vspace{-6mm}

{\footnotesize
%%%%%----------		
%%%%%---------- Abstract
\noindent{\bf Abstract}:\quad
We consider an interface with surface tension that separates a perfectly conducting inviscid fluid from a vacuum. The fluid flow is governed by the equations of ideal compressible magnetohydrodynamics (MHD), while the electric and magnetic fields in vacuum satisfy the Maxwell equations. With boundary conditions on the interface this forms a nonlinear hyperbolic problem with a characteristic free boundary. For the corresponding linearized problem we derive an energy a priori estimate in a conormal Sobolev space without assuming any stability conditions on the unperturbed flow. This verifies the stabilizing effect of surface tension because, as was shown in \cite{MT14}, a sufficiently large vacuum electric field can make the linearized problem ill-posed for the case of zero surface tension. The main ingredients in proving the energy estimate are a suitable secondary symmetrization of the Maxwell equations in vacuum and making full use of the boundary regularity enhanced from the surface tension.

\vspace{1mm}
 %%%%%----------
 %%%%%---------- Keywords

\vspace{1mm}
 %%%%%----------
 %%%%%---------- MSC2020
 %%(1980-now) Free boundary problems for PDEs

%\newpage
%%%%%----------
%%%%%---------- Introduction
%%%%%----------
\section{Introduction}\label{sec:intro}

Following \cite{CDS14,MT14}, we formulate the free interface problem when an interface separates a perfectly conducting inviscid compressible fluid (e.g., plasma) from a vacuum, and the influence of the displacement current in vacuum is taken into account.
But, unlike \cite{CDS14,MT14}, we also consider the effect of surface tension.  The fluid flow is governed by the following ideal compressible MHD equations \cite{LL84}:
\begin{align}\label{MHD}
\left\{
\begin{aligned}
&    \partial_t \rho +\div (\rho v )=0,\\
&\partial_t (\rho v )+\div (\rho v \otimes v -H \otimes H )+\nabla q =0,\\
&\partial_t H-\nabla\times(v\times H)=0,\\
&\partial_t(\rho \mathcal{E} + \tfrac{1}{2}|H|^2) +\div
\left( \rho \mathcal{E}v+pv +H\times(v\times H) \right)=0,
\end{aligned}
\right.
\end{align}
together with the divergence constraint
\begin{align} \label{divH}
 \div H=0
\end{align}
on the initial data for the unknown $ U =(q, v^{\mathsf{T}},H^{\mathsf{T}}, S)^{\mathsf{T}}\in\mathbb{R}^8$, where $\rho$ denotes density, $v=(v_1,v_2,v_3)^{\mathsf{T}}$ fluid velocity, $H =(H_1,H_2,H_3)^{\mathsf{T}} $ magnetic field, $p=p(\rho,S )$ pressure, $q =p+\frac{1}{2}|{H} |^2$ total pressure, $S$ entropy, $\mathcal{E}=\textrm{e}+\frac{1}{2}|{v}|^2$ total energy, and  $\textrm{e}=\textrm{e}(\rho,S )$ internal energy. With a state equation, $\rho=\rho(p ,S)$, and the first principle of thermodynamics, \eqref{MHD} is a closed system for $U$.

In the vacuum region, unlike the classical statement of the free interface problem in ideal compressible MHD \cite{BFKK,Goed,ST14,T10}, we do not neglect the displacement current $\varepsilon\partial_tE$ and consider the full system of Maxwell equations
\begin{align}
 \label{Maxw}
& \varepsilon\partial_th+\nabla\times E =0,\quad \varepsilon\partial_tE-\nabla\times h =0, \\
\label{Maxw-div}
& \div h=0,\quad \div E=0
\end{align}
for the magnetic field $h=(h_1,h_2,h_3)^{\mathsf{T}}$ and the electric field $E=(E_1,E_2,E_3)^{\mathsf{T}}$, where system \eqref{Maxw} is
written in a dimensionless form \cite{MT14}, and the positive constant $\varepsilon \ll 1$, being the ratio between a characteristic (average) speed of the fluid flow and the speed of light in vacuum, is a small (but {\it fixed}) parameter in the nonrelativistic setting. The divergence-free equations \eqref{Maxw-div} are constraints on the initial data for $V =(h^{\mathsf{T}},E^{\mathsf{T}})^{\mathsf{T}}\in\mathbb{R}^6$.

The problem is completed by the boundary conditions
\begin{align}
& \partial_t \varphi= v\cdot N , \label{bc1}\\
&  q- \textstyle{\frac{1}{2}}\left(|h|^2-|{E}|^2\right)=\mathfrak{s}\mathcal{H}(\varphi), \label{bc2}\\
& E\cdot{\tau_2}=\varepsilon h_3\partial_t\varphi,\quad
E\cdot{\tau_3}=-\varepsilon h_2\partial_t\varphi \label{bc3}\\
& H\cdot N=0,\quad h\cdot N =0 \label{bc4}
\end{align}
at the free interface $\Gamma (t) =\{x_1=\varphi (t,x')\}$ (with $x'=(x_2,x_3)$) separating the fluid and vacuum regions $\Omega^\pm (t)=\{\pm (x_1-\varphi (t,x'))>0\}$, where
\[
N=(1,-\p_2\varphi ,-\p_3\varphi )^{\mathsf{T}},\quad \tau_2=(\p_2\varphi , 1,0)^{\mathsf{T}},\quad \tau_3=(\p_3\varphi , 0, 1)^{\mathsf{T}},
\]
$\mathfrak{s}\geq 0$ is the
constant surface-tension coefficient,
$\mathcal{H}(\varphi)$ is twice the mean curvature of $\Gamma(t)$ defined by
\begin{align}
\nonumber	%\label{H.cal:def}
	\mathcal{H}(\varphi):= \nabla'\cdot\left( \frac{\nabla'\varphi }{\sqrt{1+|\nabla'\varphi|^2 }}\right)
	\quad \textrm{with }
	\nabla':=
	\begin{pmatrix} \p_2\\ \p_3 \end{pmatrix}.
\end{align}
Condition \eqref{bc1} means that the interface moves with the motion of the fluid whereas condition \eqref{bc2} comes from the balance of the normal stresses at the interface (for the case $\mathfrak{s}= 0$ see \cite{MT14,T12} for further discussions). Conditions \eqref{bc3} come from the jump conditions \cite{BFKK,Goed} for the electric field. At last, as in \cite{T12}, one can show that conditions \eqref{bc4} are just boundary constraints on the initial data for problem \eqref{MHD}, \eqref{Maxw}, \eqref{bc1}--\eqref{bc3}.

Note that the effect of surface tension becomes especially important in modelling the flows of liquid metals \cite{SDGX07}. The statement of problem \eqref{MHD}, \eqref{Maxw}, \eqref{bc1}--\eqref{bc3} (with constraints \eqref{divH}, \eqref{Maxw-div}, \eqref{bc4}) comes from the relativistic setting \cite{T12} of the free interface problem. As was shown in \cite{MT14}, in spite of the fact that $\varepsilon$ is a small parameter in the nonrelativistic setting, a sufficiently large vacuum electric field can make the planar interface violently unstable (for $\mathfrak{s}=0$), i.e., the corresponding constant coefficient linearized problem can be ill-posed if the unperturbed electric field in vacuum is large enough. In the classical statement \cite{BFKK,Goed,ST14,T10}, where one neglects the displacement current $\varepsilon\partial_tE$ and considers the elliptic div-curl system $\nabla\times h =0$, $\div h=0$ in the vacuum region, the influence of the vacuum electric field $E$, in contrast, is not taken into account by default. The local-in-time well-posedness of the classical free interface problem (for $\mathfrak{s}=0$) was proved in \cite{ST14}, provided that the {\it noncollinearity condition} $|H\times h|\neq 0$ holds at each point of the initial interface $\Gamma (0)$ (see also survey \cite{STT} and references therein).

The well-posedness of problem \eqref{MHD}, \eqref{Maxw}, \eqref{bc1}--\eqref{bc3} for $\mathfrak{s}= 0$ has been studied so far at the linear level. Namely, a priori estimates of the linearized problem were deduced in \cite{CDS14,MT14} for the variable coefficients satisfying the mentioned noncollinearity condition and the assumption that the normal component of the vacuum electric field is small enough. At the same time, finding a necessary and sufficient well-posedness condition for the vacuum electric field is still an open problem. This necessary and sufficient condition has been found in \cite{T20} only for the case of {\it incompressible} MHD, and under this condition satisfied for the unperturbed flow the basic $L^2$ a priori estimate was later deduced in \cite{MSTT20} for the constant coefficient linearized problem.

A simplier version of the free boundary problem \eqref{MHD}, \eqref{Maxw}, \eqref{bc1}--\eqref{bc3} with surface tension  was studied in \cite{TW22a}. This is the free boundary problem for the MHD system \eqref{MHD} in the domain $\Omega^+(t)$ with the boundary conditions \eqref{bc1}, \eqref{bc2}, \eqref{bc4} (with $\mathfrak{s}>0$) in which we formally set $h\equiv 0$ and $E\equiv 0$. The local well-posedness of this problem was established in \cite{TW22a} by a suitable modification of the Nash--Moser iteration scheme. The local well-posedness of the free boundary problem in ideal compressible MHD with surface tension and nonzero vacuum magnetic field $h$ satisfying the system $\nabla\times h =0$, $\div h=0$ was proved in \cite{TW22c} again by Nash--Moser iterations.  Since the noncollinearity condition $|H\times h|_{\Gamma (0)}\neq 0$ appearing in \cite{ST14} for $\mathfrak{s}= 0$ as well as the Taylor-type sign condition $\big(N\cdot \nabla (q-\frac{1}{2}|h|^2)\big)\big|_{\Gamma (0)}$ (which is the alternative well-posedness condition \cite{T16,TW21} for $\mathfrak{s}=0$) were not assumed in \cite{TW22c}, this verifies the stabilizing effect of surface tension on the evolution of moving vacuum
interfaces in ideal compressible MHD. Regarding the results on the well-posedness of free interface problems with surface tension in ideal incompressible MHD, we refer to \cite{GLZ21a,LL22}. Additional references as well as the references to the results for compressible and incompressible Euler equations with surface tension (see, e.g, \cite{CHS13,CS07}) can be found, for example, in \cite{STT}.

In this paper, we consider the case of positive surface tension ($\mathfrak{s}>0$) in problem \eqref{MHD}, \eqref{Maxw}, \eqref{bc1}--\eqref{bc3} and show its stabilizing role in the sense that we derive an energy a priori estimate for the corresponding linearized problem without assuming any well-posendess condition for the unperturbed flow. In particular, unlike the case $\mathfrak{s}=0$ studied in \cite{CDS14,MT14}, we do not assume that the unperturbed electric field in vacuum is small enough. The main ingredients in proving the energy estimate are a suitable secondary symmetrization \cite{T12} of the Maxwell equations in vacuum and making full use of the boundary regularity enhanced from the surface tension. Moreover, the secondary symmetrization helps to overcome the difficulty connected with the fact that the boundary is non-uniformly characteristic, i.e., characteristics are of variable multiplicity.

The plan of the rest of this paper is as follows. In Section \ref{sec:2} we reduce the free boundary problem \eqref{MHD}, \eqref{Maxw}, \eqref{bc1}--\eqref{bc3} to that in fixed domains and formulate an equivalent reduced problem with characteristics of constant multiplicity. In Section \ref{sec:3} we write down the linearized problem and state the main result for it, which is Theorem \ref{t1} about the energy a priori estimate. In Section \ref{sec:4}  we prove this estimate. At last, Section \ref{sec:5} contains concluding remarks, in particular, a brief discussion of open problems.

\section{Reduced nonlinear problem in fixed domains}
\label{sec:2}

\subsection{Symmetric hyperbolic systems of MHD and Maxwell equation}

Taking into account \eqref{divH}, we rewrite \eqref{MHD} in the nonconservative form
\begin{equation}
\label{3}
\left\{
\begin{array}{l}
{\displaystyle\frac{1}{\rho a^2}\left\{\frac{{\rm d} q}{{\rm d}t} -H
\cdot\displaystyle\frac{{\rm d} H}{{\rm d}t} \right\}+{\rm div}\,{v}=0},\qquad \rho\, {\displaystyle\frac{{\rm d}
{v}}{{\rm d} t}}-({H}\cdot\nabla ){H}+{\nabla}
q  =0 ,\\[12pt]
\displaystyle\frac{{\rm d} H}{{\rm d}t} -(H \cdot\nabla)v  -
 \frac{1}{\rho a^2}\left\{\frac{{\rm d} q}{{\rm d}t} -H
\cdot\frac{{\rm d} H}{{\rm d}t} \right\}H=0,\qquad \frac{{\rm d} S}{{\rm d} t} =0,
\end{array}
\right.
\end{equation}
where $a=a(p,S)=(\rho_p(p,S))^{-1/2}$ is the sound velocity and ${\rm d} /{\rm d} t =\partial_t+({v} \cdot{\nabla} )$. Equations \eqref{3} form the symmetric system
\begin{equation}
\label{4}
A_0(U )\partial_tU+\sum_{j=1}^3A_j(U )\partial_jU=0,
\end{equation}
with
\setlength{\arraycolsep}{4pt}
\begin{align}
\nonumber  %   \label{A0:def}
  &A_0(U):=
  \begin{pmatrix}
 \dfrac{1}{\rho a^2} & 0 & -\dfrac{1}{\rho a^2} H^{\mathsf{T}} & 0\\[2.5mm]
 0 & \rho {I}_3 & {O}_3 & 0\\[1.5mm]
 -\dfrac{1}{\rho a^2} H & {O}_3 & {I}_3+\dfrac{1}{\rho a^2}H\otimes H & 0\\[1.5mm]
 \w{0}  & \w{0}  & \w{0} & \w{1}
  \end{pmatrix},
\end{align}
and
\begin{align}
\nonumber  % \label{Ai:def}	
  &A_i(U):=
  \begin{pmatrix}
 \dfrac{v_i}{\rho a^2} & \mathbf{e}_i^{\mathsf{T}} & -\dfrac{v_i}{\rho a^2} H^{\mathsf{T}} & 0\\[2.5mm]
 \mathbf{e}_i & \rho v_i {I}_3 & -H_i {I}_3 & 0\\[1.5mm]
 -\dfrac{v_i}{\rho a^2} H & -H_i {I}_3 & v_i{I}_3+\dfrac{v_i}{\rho a^2}H\otimes H & 0\\[1.5mm]
 \w{0}  & \w{0}  & \w{0} & \w{v_i}
  \end{pmatrix}
  \ \textrm{ for $i=1,2,3$.}
\end{align}
Here and below, ${O}_m$ and ${I}_{m}$ are the zero and identity matrices of order $m$, respectively,
$\mathbf{e}_1:=(1,0,0)^{\mathsf{T}}$, $\mathbf{e}_2:=(0,1,0)^{\mathsf{T}}$, and $\mathbf{e}_3:=(0,0,1)^{\mathsf{T}}$. System \eqref{4} is hyperbolic if the matrix  $A_0$ is positive definite, i.e.,
\begin{equation}
\rho  >0,\quad \rho_p >0. \label{hypMHD}
\end{equation}
The vacuum Maxwell equations \eqref{Maxw} are also rewritten as a symmetric hyperbolic system:
\begin{equation}
\label{5}
\varepsilon\partial_tV+\sum_{j=1}^3B_j\partial_jV=0,
\end{equation}
where
\[
B_j=\begin{pmatrix}
0_3 & b_j \\
b_j^T & 0_3
\end{pmatrix},\ \textrm{ for $i=1,2,3$,}
\]
\[
b_1=\left(\begin{array}{ccc}
 0 & 0 & 0 \\
 0 & 0 & -1 \\
 0 & 1 & 0
\end{array} \right),\quad
b_2=\left(\begin{array}{ccc}
 0 & 0 & 1 \\
 0 & 0 & 0 \\
-1 & 0 & 0
\end{array} \right),\quad
b_3=\left(\begin{array}{cccccc}
 0 & -1 & 0 \\
 1 & 0 & 0 \\
 0 & 0 & 0
\end{array} \right).
\]

\subsection{An equivalent formulation in the the half-spaces}

We reformulate the free boundary problem \eqref{4}, \eqref{5}, \eqref{bc1}--\eqref{bc3} into an equivalent fixed-boundary problem by introducing the new unknowns $U_{\sharp}(t,x):=U(t,\Phi (t,x),x')$ and $V_{\sharp}(t,x):=V(t,\Phi (t,x),x')$, which are smooth in the half-spaces $\mathbb{R}^3_\pm =\{ \pm x_1>0,\ x'\in\mathbb{R}^2\}$, where
\begin{equation}\label{ch}
\Phi (t,x ):=  x_1+ \chi ( x_1)\varphi (t,x'),
\end{equation}
and $\chi\in C^{\infty}_0(-1,1)$ is the cut-off function that satisfies $\|\chi'\|_{L^{\infty}(\mathbb{R})} < 1/2$ and equals to 1 in a small neighborhood of the origin. The change of variables is non-degenerate if $\partial_1\Phi >0$. This requirement is fulfilled  if we consider solutions for which $\|\varphi \|_{L^{\infty}([0,T]\times\mathbb{R}^2)}\leq 1$. The last is true if, without loss of generality, we consider the initial data satisfying $\|\varphi_0\|_{L^{\infty}(\mathbb{R}^2)}\leq 1/2$,
and the time $T$ in our existence theorem is sufficiently small.

After the change of variables \eqref{ch} our free boundary problem \eqref{4}, \eqref{5}, \eqref{bc1}--\eqref{bc3} is reduced to the following fixed boundary problem:
\begin{subequations} \label{NP1}
	\begin{alignat}{2}
		\label{NP1a}
		&\mathbb{L}_+(U,\Phi) :=L_+(U,\Phi)U =0
		\qquad & &\quad\textrm{in  }   [0,T]\times\mathbb{R}_3^+,\\
		\label{NP1b}
		&\mathbb{L}_-(V,\Phi) :=L_-(\Phi)V =0
		& &\quad\textrm{in  }  		[0,T]\times\mathbb{R}_3^-,\\
		\label{NP1c}
		&  \mathbb{B}(U,V,\varphi)=0
		&&\quad\textrm{on  }     [0,T]\times\Gamma ,\\
		\label{NP1d}
		&U|_{t=0}=U_0,\quad V|_{t=0}=V_0,\quad \varphi|_{t=0}=\varphi_0,  &  &
	\end{alignat}
\end{subequations}
where we have dropped the subscript $``\sharp"$  for convenience, $\Gamma =\{0\}\times\mathbb{R}^2$ is the plane $x_1=0$, and
\[
L_+(U,\Phi):=A_0(U)\partial_t+\widetilde{A}_1(U,\Phi)\partial_1+ \sum_{k=2}^{3}A_k(U)\partial_k, \quad
L_-(\Phi):=\varepsilon\partial_t+\widetilde{B}_1(\Phi)\partial_1+ \sum_{k=2}^{3}B_k\partial_k,
\]
\[
\mathbb{B}(U,V,\varphi):=
\begin{pmatrix}
\partial_t \varphi- v\cdot N\\[3pt]
E\cdot{\tau_2}-\varepsilon h_3\partial_t\varphi\\[3pt]
E\cdot{\tau_3}+\varepsilon h_2\partial_t\varphi \\[3pt]
 q- \textstyle{\frac{1}{2}}|h|^2+\textstyle{\frac{1}{2}}|{E}|^2-\mathfrak{s}\mathcal{H}(\varphi)
\end{pmatrix},
\]
with
\[
\widetilde{A}_1(U,\Phi):=
	\frac{1}{\partial_1\Phi}\big(A_1(U)-\partial_t\Phi A_0(U)-\partial_2\Phi A_2(U)-\partial_3\Phi A_3(U)\big),
\]
\[
\widetilde{B}_1(\Phi) = \frac{1}{\partial_1\Phi}\big(B_1-\varepsilon\partial_t\Phi I_3-\partial_2\Phi B_2-\partial_3\Phi B_3\big).
\]
For constraints \eqref{divH}, \eqref{Maxw-div}, \eqref{bc4} one can prove the following propositions \cite{CDS14,MT14,T12}.

\begin{proposition}
Let the initial data \eqref{NP1d} satisfy
\begin{equation}
\div  \mathcal{H}=0
\label{25'}	
\end{equation}
and
\begin{equation}
H_{N}|_{x_1=0}=0,
\label{26'}
\end{equation}
where
\[
\mathcal{H}=(H_N,H_2\partial_1\Phi,H_3\partial_1\Phi),\quad
H_N=H_1-H_2\partial_2\Phi-H_3\partial_3\Phi.
\]
If problem \eqref{NP1} has a sufficiently smooth solution $(U,V,\varphi )$, then this solution satisfies \eqref{25'} and \eqref{26'} for all $t\in [0,T]$.
\label{p1}
\end{proposition}

\begin{proposition}
Let the initial data \eqref{NP1d} satisfy
\begin{equation}
{\rm div} \,\mathfrak{h}=0,\quad {\rm div}\, \mathfrak{e} =0\label{28}
\end{equation}
and
\begin{equation}
h_N|_{x_1=0}=0,\label{15.2'}
\end{equation}
where
\[
h_{N}=h_1-h_2\partial_2\Phi-h_3\partial_3\Phi,\quad \mathfrak{h}=(h_{N},h_2\partial_1\Phi,h_3\partial_1\Phi),
\]
\[
\mathfrak{e}=({E}_{N},{E}_2\partial_1\Phi,{E}_3\partial_1\Phi),\quad {E}_{N}={E}_1-{E}_2\partial_2\Phi-{E}_3\partial_3\Phi.
\]
If problem \eqref{NP1} has a sufficiently smooth solution $(U,V,\varphi )$ with the property
\begin{equation}
\partial_t\varphi \leq 0,\label{expan}
\end{equation}
then this solution satisfies  \eqref{28} and \eqref{15.2'} for all $t\in [0,T]$.
If problem \eqref{NP1} with the two additional boundary conditions
\begin{equation}
{\rm div}\, \mathfrak{h}|_{x_1=0} =0\qquad\mbox{and}\qquad
{\rm div}\, \mathfrak{e}|_{x_1=0} =0,
\label{add.bound}
\end{equation}
has a sufficiently smooth solution $(U,V,\varphi )$ with the property
\begin{equation}
\partial_t\varphi > 0,\label{shrink}
\end{equation}
then this solution again satisfies \eqref{28} and \eqref{15.2'} for all $t\in [0,T]$.
\label{p2}
\end{proposition}

\subsection{Correct number of the boundary conditions}

As in \cite{T10}, for the MHD system \eqref{NP1a}, the boundary matrix $\widetilde{A}_1(U,\Phi )|_{x_1=0}$ on the boundary $\Gamma$ has one positive and one negative eigenvalue and the others are zero. This means that the boundary $\Gamma$ is {\it characteristic}, and according to the number of incoming/outgoing characteristics, the hyperbolic system \eqref{NP1a} in the half-space $\mathbb{R}^3_+$ requires one boundary condition on $\Gamma$. At the same time, for the Maxwell system \eqref{NP1b}, the boundary matrix $B=B(\varphi):=\widetilde{B}_1({\Phi})|_{x_1=0}$  has the eigenvalues
\begin{align*}
&
\lambda_{1,2}(B)= -\varepsilon\partial_t\varphi -
\sqrt{1 +(\partial_2\varphi)^2+(\partial_3\varphi)^2},\qquad \lambda_{3,4}(B) =-\varepsilon\partial_t\varphi , \\
& \lambda_{5,6}(B)= -\varepsilon\partial_t\varphi +
\sqrt{1 +(\partial_2\varphi)^2+(\partial_3\varphi)^2}.
\end{align*}
If \eqref{expan} holds, then the matrix $B$ has two negative eigenvalues (recall that $\varepsilon \ll 1$). This implies that the hyperbolic system \eqref{NP1b} in $\mathbb{R}^3_-$ requires two boundary conditions on $\Gamma$. That is, under assumption \eqref{expan} we have a correct number of boundary conditions in \eqref{NP1c} because the first condition in \eqref{NP1c}  is needed for determining the function $\varphi$.

If the plasma region is shrinking, i.e., \eqref{shrink} holds, then the matrix $B$ has four negative eigenvalues. Hence, the correct number of boundary conditions is six, and problem \eqref{NP1} is missing two boundary conditions. However, if, following \cite{T12}, we supplement \eqref{NP1c} with the additional boundary conditions \eqref{add.bound}, which enables one to prove \eqref{28}, we have a correct number of boundary conditions also for case \eqref{shrink}. Completing our problem with
\begin{equation}\label{addbc}
 {\rm div}\, \mathfrak{h} =0,\quad {\rm div}\, \mathfrak{e} =0\qquad\mbox{on } [0,T]\times\Gamma\cap\{\p_t\varphi >0\},
\end{equation}
we come to problem \eqref{NP1}, \eqref{addbc}, which is well formulated according to the number of boundary conditions.

\subsection{Equivalent problem with characteristics of constant multiplicity}

Since we prescribe a different number of boundary conditions on different portions of the boundary $\Gamma$, the boundary is non-uniformly characteristic. For overcoming this difficulty it was proposed in \cite{CDS14} to introduce such a new unknown in the vacuum region that
the boundary becomes characteristic of constant multiplicity.  Here, as in \cite{T24}, we prefer to follow an alternative way consisting in the application of a secondary symmetrization to the vacuum Maxwell equations. This symmetrization  proposed in \cite{T12} will, first, enable us to reduce our nonlinear \eqref{NP1}, \eqref{addbc}  to that with characteristics of constant multiplicity and, second, derive, as in \cite{CDS14,MT14,T12}, an energy a priori estimate for the linearized problem.

The {\it secondary symmetrization} \cite{T12} of the symmetric system \eqref{5} reads
\begin{equation}
\varepsilon\mathcal{B}_0\partial_tV +\sum_{j=1}^3\mathcal{B}_j\partial_jV=0,
\label{20}
\end{equation}
where
\[
\mathcal{B}_0=\left(\begin{array}{cccccc}
1 & 0 & 0& 0 & \nu_3 & -\nu_2 \\
0 & 1 & 0& -\nu_3 & 0 & \nu_1 \\
0 & 0 & 1& \nu_2 & -\nu_1 & 0 \\
0 & -\nu_3 & \nu_2& 1 & 0 & 0 \\
\nu_3 & 0 & -\nu_1& 0 & 1 & 0 \\
-\nu_2 & \nu_1 & 0& 0 & 0 & 1
\end{array} \right),\quad
\mathcal{B}_1=
\left(\begin{array}{cccccc}
\nu_1 & \nu_2 & \nu_3& 0 & 0 & 0 \\
\nu_2 & -\nu_1 & 0& 0 & 0 & -1 \\
\nu_3 & 0 & -\nu_1& 0 & 1 & 0 \\
0 & 0 & 0& \nu_1 & \nu_2 & \nu_3 \\
0 & 0 & 1& \nu_2 & -\nu_1 & 0 \\
0 & -1 & 0& \nu_3 & 0 & -\nu_1
\end{array} \right),
\]
\[
\mathcal{B}_2=
\left(\begin{array}{cccccc}
-\nu_2 & \nu_1 & 0& 0 & 0 & 1 \\
\nu_1 & \nu_2 & \nu_3& 0 & 0 & 0 \\
0 & \nu_3 & -\nu_2& -1 & 0 & 0 \\
0 & 0 & -1& -\nu_2 & \nu_1 & 0 \\
0 & 0 & 0& \nu_1 & \nu_2 & \nu_3 \\
1 & 0 & 0& 0 & \nu_3 & -\nu_2
\end{array} \right),\quad
\mathcal{B}_3=
\left(\begin{array}{cccccc}
-\nu_3 & 0 & \nu_1& 0 & -1 & 0 \\
0 & -\nu_3 & \nu_2& 1 & 0 & 0 \\
\nu_1 & \nu_2 & \nu_3& 0 & 0 & 0 \\
0 & 1 & 0& -\nu_3 & 0 & \nu_1 \\
-1 & 0 & 0& 0 & -\nu_3 & \nu_2 \\
0 & 0 & 0& \nu_1 & \nu_2 & \nu_3
\end{array} \right),
\]
and  $\nu_1$, $\nu_2$ and $\nu_3$ are arbitrary functions of  $(t,x)$.
System \eqref{20} is equivalent to \eqref{5} and it is again hyperbolic if
$\mathcal{B}_0>0$, i.e.,
\begin{equation}
|\nu |<1,\label{21'}
\end{equation}
with the vector-function $\nu = (\nu_1, \nu_2, \nu_3)$. Using the divergence constraint \eqref{28} and omitting technical details (see \cite{CDS14,MT14,ST13}), for system \eqref{NP1b} we obtain the following counterpart of the secondary symmetrization \eqref{20}:
\begin{equation}
\varepsilon\mathcal{B}_0(\nu )\partial_t{V}+\widetilde{\mathcal{B}}_1(\nu, \Phi )\partial_1{V}+\sum_{k=2}^3\mathcal{B}_k(\nu )\partial_k{V}=0 \qquad \mbox{in}\ [0,T]\times\mathbb{R}_3^-,\label{105v}
\end{equation}
where
\[
\widetilde{\mathcal{B}}_1(\nu ,\Phi )= \frac{1}{\partial_1{\Phi}}\left(\mathcal{B}_1(\nu )-\varepsilon \partial_t\Phi \mathcal{B}_0(\nu ) -\partial_2{\Phi}\mathcal{B}_2(\nu ) -
\partial_3{\Phi}\mathcal{B}_3(\nu )\right).
\]
We now make the following choice of  $\nu (t,x )$ suggested by the choice made in \cite{CDS14,MT14} for the linearized problem with zero surface tension ($\mathfrak{s}= 0$):
\begin{equation}
\label{ch1}
\nu = \varepsilon v^-,
\end{equation}
where $v^-=v^-(t,x):=v(t,-x_1,x')$.
Since $\varepsilon$ is a small parameter, the hyperbolicity condition \eqref{21'} holds for choice \eqref{ch1}. Moreover, one can show that smooth solutions of \eqref{105v} satisfy the divergence constraints \eqref{28} for all $t\in [0,T]$ if they were true at $t=0$ (see \cite{CDS14,ST13} for the proof). This implies the equivalence of problem \eqref{NP1}, \eqref{addbc} and problem \eqref{NP1a}, \eqref{NP1c}, \eqref{NP1d}, \eqref{105v}.

With the notations
\begin{equation}
\mathfrak{B}_k(v^-):=\mathcal{B}_k(\varepsilon v^- ),\quad k=0,2,3,\qquad \mathfrak{B}_1(v^-,\Phi):=
\widetilde{\mathcal{B}}_1(\varepsilon v^- ,\Phi ),\label{bm}
\end{equation}
system \eqref{105v} is rewritten as
\begin{equation}
\mathbb{L}_-(v^-,V,\Phi) :=L_-(v^-,\Phi)V=0
 \qquad \mbox{in}\ [0,T]\times\mathbb{R}_3^-,\label{NP1b'}
\end{equation}
where
\[
L_-(v^-,\Phi):=
\varepsilon\mathfrak{B}_0(v^-)\partial_t+\mathfrak{B}_1(v^-,\Phi)\partial_1+\mathfrak{B}_2(v^-)\partial_2+
\mathfrak{B}_3(v^-)\partial_3.
\]
Taking the first boundary condition in \eqref{NP1c} into account, we compute the eigenvalues of the boundary matrix $\mathcal{B}=\mathcal{B}(v|_{x_1=0},\varphi):=\mathfrak{B}_1(v^-,\Phi)|_{x_1=0}$:
\begin{align*}
&
\lambda_{1,2}(\mathcal{B})= - \sqrt{1 +(\partial_2\varphi)^2+(\partial_3\varphi)^2} +\mathcal{O} (\varepsilon ),\quad \lambda_{3,4}(\mathcal{B}) =0 ,\\
& \lambda_{5,6}(\mathcal{B})=  \sqrt{1 +(\partial_2\varphi)^2+(\partial_3\varphi)^2} +\mathcal{O} (\varepsilon ).
\end{align*}
As we can see, the hyperbolic system \eqref{NP1b'}  requires two boundary conditions on $\Gamma$. This means that problem \eqref{NP1a}, \eqref{NP1c}, \eqref{NP1d}, \eqref{NP1b'} has a correct number of boundary conditions in \eqref{NP1c} regardless of the sign of $\partial_t\varphi$.
From now on we will consider the initial boundary value problem \eqref{NP1a}, \eqref{NP1c}, \eqref{NP1d}, \eqref{NP1b'} for which the boundary is characteristic of constant multiplicity.

\section{Linearized problem and main result}
\label{sec:3}

\subsection{Basic state}

Let
\begin{equation}
(\mathring{U}(t,x ),\mathring{V}(t,x ),\mathring{\varphi}(t,{x}'))
\label{37}
\end{equation}
be a given sufficiently smooth vector-function with $\mathring{U}=(\mathring{q},\mathring{v}^{\mathsf{T}},\mathring{H}^{\mathsf{T}},\mathring{S})^{\mathsf{T}}$, $\mathring{V}=(\mathring{h}^{\mathsf{T}},\mathring{E})^{\mathsf{T}}$, and
\[
\|\mathring{U}\|_{W^3_{\infty}(\Omega_T^+)}+
\|\mathring{V}\|_{W^3_{\infty}(\Omega_T^-)}+
\|\mathring{\varphi}\|_{W^4_{\infty}(\Gamma_T)} \leq K,
\]
where $K>0$ is a constant,
\[
\Omega_T^\pm:= (-\infty, T]\times\mathbb{R}^3_\pm,\quad
\Gamma_T:=(-\infty ,T]\times\Gamma,
\]
and below all the ``ring'' values like $\mathring{U}$ will be related to the {\it basic state} \eqref{37}. Following \cite{CDS14,MT14,T12}, we also assume that the basic state \eqref{37} satisfies the hyperbolicity conditions \eqref{hypMHD}, the first  three boundary conditions in \eqref{NP1c}, the equations for $H$ and $h$ contained in \eqref{NP1a} and \eqref{NP1b}, constraints \eqref{25'}--\eqref{15.2'} at $t=0$, and the inequality $\|\mathring{\varphi} \|_{L^{\infty}([0,T]\times\mathbb{R}^2)}\leq 1$.

\subsection{Linearized problem}

The linearized operators for problem \eqref{NP1a}, \eqref{NP1c}, \eqref{NP1d}, \eqref{NP1b'} are defined as follows:
\[
\begin{aligned}
&\mathbb{L}_+'\big(\mathring{U},\mathring{\Phi}\big)(U,\Phi)
:=\left.\frac{d}{d\theta}
\mathbb{L}_+\big(\mathring{U}+{\theta}U, \mathring{\Phi}+{\theta}\Phi\big)\right|_{\theta=0},\\[6pt]
&\mathbb{L}_-'\big(\mathring{W},\mathring{\Phi}\big)(W,\Phi)
:=\left.\frac{d}{d\theta}
\mathbb{L}_{-}\big(\mathring{W}+{\theta}W, \mathring{\Phi}+{\theta}\Phi\big)\right|_{\theta=0},\\[6pt]
&\mathbb{B}'\big(\mathring{U},\mathring{V},\mathring{\varphi}\big)(U,V,\varphi)
:=\left.\frac{d}{d\theta}
\mathbb{B}(\mathring{U}+{\theta}U,\mathring{V}+{\theta}V,
\mathring{\varphi}+{\theta}\varphi)\right|_{\theta=0},
\end{aligned}
\]
where $\mathring{\Phi}(t,x):=x_1+\chi(x_1)\mathring{\varphi}(t,x') $, $\mathring{W}:=(\mathring{V}^{\mathsf{T}},\mathring{v}^-)^{\mathsf{T}}$, ${W}:=({V}^{\mathsf{T}},{v}^-)^{\mathsf{T}}$ and $\mathring{v}^-=\mathring{v}^-(t,x):=\mathring{v}(t,-x_1,x')$. We easily compute them. In particular, the linearized interior equations read:
\begin{equation}
\begin{aligned}
& \mathbb{L}_+'\big(\mathring{U},\mathring{\Phi}\big)(U,\Phi)
 = \mathbb{L}'_{e+}\big(\mathring{U},\mathring{\Phi}\big)U
-\frac{L_+(\mathring{U}, \mathring{\Phi})\Phi}{\p_1 \mathring{\Phi}}\,\p_1\mathring{U},\\[6pt]
& \mathbb{L}_-'\big(\mathring{W},\mathring{\Phi}\big)(W,\Phi)
 = \mathbb{L}'_{e-}\big(\mathring{W},\mathring{\Phi}\big)W
-\frac{L_-(\mathring{v}^-, \mathring{\Phi})\Phi}{\p_1 \mathring{\Phi}}\,\p_1\mathring{V},
\end{aligned}
\label{op}
\end{equation}
where
\[
\mathbb{L}'_{e+}\big(\mathring{U},\mathring{\Phi}\big)U:=L_+\big(\mathring{U},\mathring{\Phi}\big)U+
\mathcal{C}_{+}( \mathring{U},\mathring{\Phi})U,\qquad  \mathbb{L}'_{e-}\big(\mathring{W},\mathring{\Phi}\big)W:=L_-\big(\mathring{v}^-,\mathring{\Phi}\big)V+
\mathcal{C}_{-}( \mathring{V},\mathring{\Phi})v^-,
\]
and the concrete form of the matrices $\mathcal{C}_{\pm}$ is of no interest (see \cite{CDS14,MT14,T12}).

The differential operators $\mathbb{L}_\pm'$ are first-order operators in $\Phi$, i.e., the linearized interior equations contain derivatives of the interface perturbation. For getting standard linear hyperbolic systems we first pass to the Alinhac's good unknowns \cite{A89}
\begin{align} \label{good}
	\dot{U}=(\dot{q}, \dot{v}^{\mathsf{T}},\dot{H}^{\mathsf{T}}, \dot{S})^{\mathsf{T}}:=U-\frac{\Psi}{\partial_1 \mathring{\Phi}}\partial_1\mathring{U},
	\qquad \dot{V}=(\dot{h}^{\mathsf{T}},\dot{E})^{\mathsf{T}}:=V-\frac{\Psi}{\partial_1 \mathring{\Phi}}\partial_1\mathring{V},
\end{align}
with $\Psi(t,x):=\chi(x_1)\varphi(t,x')$. In terms of \eqref{good} the operators in \eqref{op} are rewritten as
\begin{equation}
\begin{aligned}
& \mathbb{L}_+'\big(\mathring{U},\mathring{\Phi}\big)(U,\Phi)
 = \mathbb{L}'_{e+}\big(\mathring{U},\mathring{\Phi}\big)\dot{U}
+\frac{\Psi}{\partial_1\mathring{\Phi}}\,
\partial_1\mathbb{L}_{+}(\mathring{U},\mathring{\Phi} ),\\[6pt]
& \mathbb{L}_-'\big(\mathring{W},\mathring{\Phi}\big)(W,\Phi)
 = \mathbb{L}'_{e-}\big(\mathring{W},\mathring{\Phi}\big)\dot{W}+\frac{\Psi}{\partial_1\mathring{\Phi}}\,
\partial_1\mathbb{L}_{-}(\mathring{W}, \mathring{\Phi}),
\end{aligned}
\label{op'}
\end{equation}
where $\dot{W}:=(\dot{V}^{\mathsf{T}},\dot{v}^-)^{\mathsf{T}}$ and $\dot{v}^-=\dot{v}^-(t,x):=\dot{v}(t,-x_1,x')$. Then, we drop the zero-order terms in $\Psi$  in \eqref{op'}, which in the future nonlinear analysis are considered as error terms at each Nash--Moser iteration step. This gives us the following final form of our linearized problem for $(\dot{U},\dot{V},\varphi )$:
\begin{subequations} \label{ELP1}
	\begin{alignat}{2}
		&L_+\big(\mathring{U},\mathring{\Phi}\big)\dot{U}+
\mathcal{C}_{+}( \mathring{U},\mathring{\Phi})\dot{U} =f
		&\quad  &\textrm{in } \Omega^+_T,
		\label{ELP1a}\\
		&L_-\big(\mathring{v}^-,\mathring{\Phi}\big)\dot{V}+
\mathcal{C}_{-}( \mathring{V},\mathring{\Phi})\dot{v}^-=0
		&\quad  &\textrm{in } \Omega^-_T,
		\label{ELP1b}\\
		\label{ELP1c}   & 		\mathbb{B}'\big(\mathring{U},\mathring{V},\mathring{\varphi}\big)(\dot{U},\dot{V},\varphi)
		=0\quad
		&&\textrm{on } \Gamma_{T},\\
		&(\dot{U},\dot{V},\varphi)\big|_{t<0}=0,
		&&
		\label{ELP1d}
	\end{alignat}
\end{subequations}
where
\[
\mathbb{B}'\big(\mathring{U},\mathring{V},\mathring{\varphi}\big)(\dot{U},\dot{V},\varphi):=
\begin{pmatrix}
\big(\p_t + \mathring{v}'\cdot\nabla'  -\p_1 (\mathring{v}\cdot\mathring{N})\big) \varphi-\dot{v}\cdot\mathring{N} \\[3pt]
\dot{E}\cdot\mathring{\tau}_2 -\varepsilon (\p_t\mathring{\varphi}) \dot{h}_3-\varepsilon\p_t(\mathring{h}_3\varphi ) +\p_2(\mathring{E}_1\varphi)\\[3pt]
\dot{E}\cdot\mathring{\tau}_3 +\varepsilon (\p_t\mathring{\varphi}) \dot{h}_2+\varepsilon\p_t(\mathring{h}_2\varphi) +\p_3(\mathring{E}_1\varphi)\\[3pt]
 \dot{q}- \mathring{h}\cdot\dot{h} +\mathring{E}\cdot\dot{E}+[\partial_1\mathring{q}] \varphi -\mathfrak{s}\nabla'\cdot (\mathring{B}\nabla'\varphi)
\end{pmatrix},
\]
\[
\mathring{v}'=(\mathring{v}_2,\mathring{v}_3)^{\mathsf{T}},\quad
\nabla'=(\p_2,\p_3)^{\mathsf{T}},\quad \mathring{N}=(1,-\p_2\mathring{\varphi} ,-\p_3\mathring{\varphi} )^{\mathsf{T}},
\]
\[
\mathring{\tau}_2=(\p_2\mathring{\varphi} , 1,0)^{\mathsf{T}},\quad \mathring{\tau}_3=(\p_3\mathring{\varphi} , 0, 1)^{\mathsf{T}},\quad
[\partial_1\mathring{q}]=(\partial_1\mathring{q})|_{\Gamma}-(\mathring{h}\cdot
\partial_1\mathring{h})|_{\Gamma}+(\mathring{E}\cdot\partial_1\mathring{E})|_{\Gamma},
\]
and $\mathring{B}$ is the positive definite matrix defined by (see \cite{TW22c})
\begin{align}
\label{B.ring:def}
  \mathring{B}:=
\frac{I_2}{|\mathring{N}|}
- \frac{\nabla'\mathring{\varphi}\otimes \nabla'\mathring{\varphi}}{|\mathring{N}|^3}.
\end{align}
The assumption that the basic state \eqref{37} satisfies the equation $h$ contained in \eqref{NP1b} was used while writing down the second and third boundary conditions in \eqref{ELP1c} (see see \cite{CDS14,MT14,T12}). We assume that the given source term $f$ vanishes in the past and consider the case of zero initial data, which is the usual assumption. The case of nonzero initial data is postponed to the nonlinear analysis (construction of a so-called approximate solution; see, e.g., \cite{ST14}).

Moreover, we consider the homogeneous equations in system \eqref{ELP1b} (with zero source terms) and the homogeneous boundary conditions \eqref{ELP1c} because, following  \cite{CDS14,MT14,T12}, the linearized problem with inhomogeneous vacuum equations and inhomogeneous boundary conditions can be reduced to problem \eqref{ELP1}. Note that the process of reduction of the linearized problem to that with homogeneous vacuum equations and homogeneous boundary conditions described in \cite{CDS14,T12} for $\mathfrak{s}=0$ is the same for our case when $\mathfrak{s}\neq 0$ (see \cite{CDS14,T12} for more details). It is worth noting that this process is organized so that the solutions of the reduced problem \eqref{ELP1} automatically satisfy the following  linear versions of constraints \eqref{25'}--\eqref{15.2'}:
\begin{alignat}{2}
		&\div  \dot{\mathcal{H}}=0
		& \quad &\textrm{in } \Omega^+_T,
		\label{div+}\\
		&\div\dot{\mathfrak{h}}=0,\quad \div\dot{\mathfrak{e}}=0
		&\quad  &\textrm{in } \Omega^-_T,
		\label{div-}\\
		\label{H_N}   & 		
\dot{H}_N=\mathring{H}_2\partial_2\varphi +\mathring{H}_3\partial_3\varphi - \varphi\,\partial_1(\mathring{H}\cdot\mathring{N})\quad
		& &\textrm{on } \Gamma_{T},\\
& \dot{h}_N=\mathring{h}_2\partial_2\varphi +\mathring{h}_3\partial_3\varphi - \varphi\,\partial_1(\mathring{h}\cdot\mathring{N})
		& &\textrm{on } \Gamma_{T},
		\label{h_N}
	\end{alignat}
where
\[
\dot{\mathcal{H}}=(\dot{H}_N,\dot{H}_2\partial_1\mathring{\Phi},
\dot{H}_3\partial_1\mathring{\Phi})^{\mathsf{T}}, \quad \dot{H}_N=\dot{H}_1-\dot{H}_2\partial_2\mathring{\Phi}-\dot{H}_3\partial_3\mathring{\Phi},
\quad \dot{H}_N|_{\Gamma}= (\dot{H}\cdot\mathring{N})|_{\Gamma},
\]
\[
\dot{\mathfrak{h}}=(\dot{h}_{N},
\dot{h}_2\partial_1\mathring{\Phi},
\dot{h}_3\partial_1\mathring{\Phi})^{\mathsf{T}},\quad
\dot{h}_N=\dot{h}_1-\dot{h}_2\partial_2\mathring{\Phi}-\dot{h}_3\partial_3\mathring{\Phi},
\quad \dot{h}_N|_{\Gamma}= (\dot{h}\cdot\mathring{N})|_{\Gamma},
\]
\[
\dot{\mathfrak{e}}=(\dot{E}_{N},
\dot{E}_2\partial_1\mathring{\Phi},
\dot{E}_3\partial_1\mathring{\Phi})^{\mathsf{T}},\quad
\dot{E}_N=\dot{E}_1-\dot{E}_2\partial_2\mathring{\Phi}-\dot{E}_3\partial_3\mathring{\Phi}.
\]

\subsection{Conormal Sobolev spaces}

Since the boundary $\Gamma$ is characteristic, we have a natural loss of control on derivatives of the unknowns in the normal direction. At the same time, for the vacuum unknown $\dot{V}$ we can compensate this loss thanks the structure of the boundary matrix and the divergence-free equations \eqref{div-}. We now provide the natural functional setting for the plasma unknown $\dot{U}$. We define the conormal derivative $\mathrm{D}_{\rm tan}^{\alpha}$ by
\begin{align}
\label{D*:def}
\mathrm{D}_{\rm tan}^{\alpha}:=\p_t^{\alpha_0} (\sigma \p_1)^{\alpha_1}\p_2^{\alpha_2}  \p_3^{\alpha_{3}}\qquad\textrm{for }  \alpha:=(\alpha_0,\ldots,\alpha_{3})\in\mathbb{N}^{4},	
\end{align}
where $\sigma (x_1)\in C^{\infty}(\mathbb{R}_+)$ is
a monotone increasing function such that $\sigma (x_1)=x_1$ in a neighborhood of
the origin and $\sigma (x_1)=1$ for $x_1$ large enough.
For $m\in\mathbb{N}$, the conormal Sobolev space $H_{\rm tan}^{m}(\Omega^+_T)$ (see \cite{MS}) is defined as
\begin{align*}
	&H_{\rm tan}^m(\Omega^+_T):=
	\{ u\in L^2(\Omega^+_T):\, \mathrm{D}_{\rm tan}^{\alpha} u\in L^2(\Omega^+_T)
	\textrm{ for } | \alpha |\leq m  \},
\end{align*}
and equipped with the norm
\begin{align}  \nonumber
	{\|}u{\|}^2_{H^m_{\rm tan}(\Omega^+_T)}:=
	\sum_{| \alpha|\leq m} \|\mathrm{D}_{\rm tan}^{\alpha} u\|_{L^2(\Omega^+_T)}^2.
\end{align}
Below  we will also use the norms
\[
\nt u(t)\nt^2_{H^m_{\rm tan}(\mathbb{R}^3_+)}:=\sum_{| \alpha|\leq m} \|\mathrm{D}_{\rm tan}^{\alpha} u(t)\|_{L^2(\mathbb{R}^3_+)}^2\quad\mbox{and}\quad
\nt u(t)\nt^2_{H^m(\mathbb{R}^3_-)}:=\sum_{| \alpha|\leq m} \|\mathrm{D}^{\alpha} u(t)\|_{L^2(\mathbb{R}^3_-)}^2,
\]
where $\mathrm{D}^{\alpha}:=\p_t^{\alpha_0} \p_1^{\alpha_1}\p_2^{\alpha_2}  \p_3^{\alpha_{3}}$.

\subsection{Main result}

Hereafter, we use $A\lesssim_{a_1,\ldots,a_m} B$ to denote that $A \leq C(a_1,\ldots,a_m)B$ for given parameters $a_1,\ldots,a_m$, where we denote by $C$ some universal positive constant
and by $C(\cdot)$ some positive constant depending on the quantities listed in the parenthesis.

We are now in the position to state the main result of this paper.

\begin{theorem}
Let the basic state \eqref{37} satisfies the assumptions formulated above and the source term $f \in H^1_{\rm tan}(\Omega_T^+)$ vanishes in the past. Let
 problem \eqref{ELP1}  has a solution $(\dot{U},\dot{V},\varphi)\in H^1_{\rm tan}(\Omega_T^+)\times H^1(\Omega_T^-)\times H^{1}(\Gamma_T)$, with $\nabla'\varphi \in H^1(\Gamma_T)$. Then this solution obeys
the a priori estimate
\begin{equation}
\|\dot{U}\|_{H^{1}_{\rm tan}(\Omega_T^+)}+\|\dot{V}\|_{H^{1}(\Omega_T^-)}+\|(\varphi , \nabla'\varphi )\|_{H^{1}(\Gamma_T)}\lesssim_{K,T}
\|f\|_{H^{1}_{\rm tan}(\Omega_T^+)}.
\label{54}
\end{equation}
\label{t1}
\end{theorem}

\section{Energy a priori estimate}
\label{sec:4}

\subsection{Energy inequality in $L^2$}

By standard arguments of the energy method applied to the symmetric hyperbolic systems \eqref{ELP1a} and \eqref{ELP1b}, we obtain
\begin{equation}
I(t)+ \int_{\Gamma_t}\mathcal{Q}\lesssim
\| f\|^2_{L_2(\Omega_T^+)} +\|\dot{U}\|_{L^2 (\Omega_t^+)}+\|\dot{V}\|_{L^2(\Omega_t^-)},
\label{107}
\end{equation}
where
\[
I(t)= \int_{\mathbb{R}^3_+}A_0(\mathring{U})\dot{U}\cdot \dot{U} +\int_{\mathbb{R}^3_-}\mathfrak{B}_0(\mathring{v})\dot{V}\cdot\dot{V},\quad
\mathcal{Q}  =-\big(\widetilde{A}_1(\mathring{U},\mathring{\Phi})
\dot{U}\cdot \dot{U}\big)\big|_{\Gamma}+
\frac{1}{\varepsilon}\big(\mathfrak{B}_1(\mathring{v}^-,\mathring{\Phi})
\dot{V}\cdot\dot{V}\big)\big|_{\Gamma}\,,
\]
in particular, we have
\begin{equation}
\big(\widetilde{A}_1(\mathring{U},\mathring{\Phi})\dot{U}\cdot \dot{U}\big)\big|_{\Gamma}=2\dot{q}\dot{v}_N|_{\Gamma},
\label{AAA}
\end{equation}
with $\dot{v}_N=\dot{v}_1-\dot{v}_2\partial_2\mathring{\Phi}-\dot{v}_3\partial_3\mathring{\Phi}$ (clearly, $\dot{v}_N|_{\Gamma}=(\dot{v}\cdot\mathring{N})|_{\Gamma}$).
Thanks to the choice in \eqref{ch1}, by using the boundary conditions it was shown in \cite{CDS14,MT14,T12} that the quadratic form $\mathcal{Q}$ for $\mathfrak{s}=0$ can be reduced to the form
\begin{equation}
\mathcal{Q} = \p_t\big(\mathring{\mu}\varphi \dot{E}_N\big)+\p_2\big(\mathring{\mu}\varphi (\dot{E}_2\p_t\mathring{\varphi}-\dot{h}\cdot\mathring{\tau}_3)\big)+\p_3\big(\mathring{\mu}\varphi (\dot{E}_3\p_t\mathring{\varphi}+\dot{h}\cdot\mathring{\tau}_2)\big) +\mathcal{L}\qquad\mbox{on }\Gamma ,
\label{QQ}
\end{equation}
where $\mathring{\mu}=2\big(\mathring{E}_1+\varepsilon\mathring{v}_2\mathring{h}_3-
\varepsilon\mathring{v}_3\mathring{h}_2\big)$ and
$\mathcal{L}$ is the sum of terms like ${\rm coeff}\,\dot{q}\varphi$,
${\rm coeff}\,\dot{v}_N\varphi$, ${\rm coeff}\,\dot{h}_i\varphi$, ${\rm coeff}\,\dot{E}_i\varphi$, ${\rm coeff}\,\varphi^2$. Here and below coeff is a
generic coefficient, which depends on the basic state \eqref{37}, whose exact form is of no interest and it may change from line to line. For $\mathfrak{s}>0$, we have the additional term
\[
\mathcal{S}=
-\mathfrak{s}\big(\nabla'\cdot (\mathring{B}\nabla'\varphi)\big) \big(\p_t + \mathring{v}'\cdot\nabla'  -\p_1 (\mathring{v}\cdot\mathring{N})\big) \varphi
\]
in the right-hand side of \eqref{QQ}, which comes from \eqref{AAA} if we use the first and last boundary conditions in \eqref{ELP1c}:
\begin{equation}
\mathcal{Q} = \mathcal{S}+\p_t\big(\mathring{\mu}\varphi \dot{E}_N\big)+\p_2\big(\mathring{\mu}\varphi (E_2\p_t\mathring{\varphi}-\dot{h}\cdot\mathring{\tau}_3)\big)+\p_3\big(\mathring{\mu}\varphi (E_3\p_t\mathring{\varphi}+\dot{h}\cdot\mathring{\tau}_2)\big) +\mathcal{L}\qquad\mbox{on }\Gamma .
\label{QQQ}
\end{equation}

Referring to \cite[(2.20)]{TW22a}, we rewrite the term $\mathcal{S}$ (responsible for the surface tension) as follows:
\begin{equation}
\mathcal{S}=\mathfrak{s}\p_t\big(\mathring{B}\nabla'\varphi\cdot\nabla'\varphi
\big) +\p_2\{\ldots\} +\p_3\{\ldots\}+Q_0(\varphi ,\nabla'\varphi  ),
\label{QQQQ}
\end{equation}
where the concrete form of the expressions in curly braces and the quadratic form $Q_0$ of the variable $(\varphi ,\nabla'\varphi )$ is of no interest. Using \eqref{B.ring:def},
from \eqref{107}, \eqref{QQQ}, and \eqref{QQQQ} we deduce the energy inequality
\begin{multline}
I(t)+\int_{\Gamma}\bigg\{\frac{\mathfrak{s}|\nabla'\varphi |^2}{|\mathring{N}|^3} + \mathring{\mu}\varphi \dot{E}_N\bigg\}  + \int_{\Gamma_t}\mathcal{L} \\
\lesssim_K \|f\|_{L^2(\Omega_T^+)} +\|\dot{U}\|_{L^2 (\Omega_t^+)}+\|\dot{V}\|_{L^2(\Omega_t^-)}+\|(\varphi ,\nabla'\varphi )\|^2_{L^2(\Gamma_t )} .\qquad
\label{enegineq}
\end{multline}
As we can see, we are not able to close the a priori estimate in $L^2$.

\subsection{Preparatory estimates}

For closing the estimate in $H^1$ (in the sense of \eqref{54}) we need some preparatory estimates. Using the special structure of the boundary matrix $\widetilde{A}_1(\mathring{U},\mathring{\Phi})$ (see \eqref{AAA}), from system \eqref{ELP1a} we deduce the estimate
\begin{equation}
\|(\partial_1\dot{q},\partial_1\dot{v}_N) (t)\|^2_{L^2(\mathbb{R}^3_+)}\lesssim_K \|f(t)\|^2_{L^2(\mathbb{R}^3_+)}+\nt \dot{U}(t)\nt^2_{H^1_{\rm tan}(\mathbb{R}^3_+)},
\label{qvn}
\end{equation}
which implies the estimate
\begin{equation}
\|(\dot{q},\dot{v}_N) (t)\|^2_{L^2(\Gamma )}\lesssim_K \|f(t)\|^2_{L^2(\mathbb{R}^3_+)}+\nt \dot{U}(t)\nt^2_{H^1_{\rm tan}(\mathbb{R}^3_+)}
\label{qvntr}
\end{equation}
for the traces $\dot{q}|_{\Gamma}$ and $\dot{v}_N|_{\Gamma}$. Since $\sigma\p_1\dot{U}|_{\Gamma}=0$, we do not need to use the boundary conditions for estimating the weighted derivative $\sigma\p_1\dot{U}$, and we easily obtain
\begin{equation}
\|\sigma\p_1\dot{U}(t)\|_{L^2(\Omega_t^+)}\lesssim_K
\|(f,\dot{U})\|_{H^1_{\rm tan}(\Omega_t^+)}.
\label{wder}
\end{equation}

Thanks to the special structure of the boundary matrix $\mathfrak{B}_1(\mathring{v}^-,\mathring{\Phi})$ (see \cite{CDS14,MT14,T12}) and the divergences \eqref{div-} we have the full control on the normal derivatives for $\dot{V}$ in the sense that
\begin{equation}
\|\partial_1\dot{V} (t)\|_{L^2(\mathbb{R}^3_-)}\lesssim_K
\nt \dot{V}(t)\nt_{H^1(\mathbb{R}^3_-)}.
\label{d1V}
\end{equation}
It follows from the trace theorem that
\begin{equation}
\|\dot{V} (t)\|_{L^2(\Gamma )}\lesssim_K \nt \dot{V}(t)\nt_{H^1(\mathbb{R}^3_-)}.
\label{Vtr}
\end{equation}

Below we will use the elementary inequalities
\begin{equation}\label{el1}
\|\dot{U}(t)\|_{L^2(\mathbb{R}^3_+)}\leq \|(\dot{U},\p_t\dot{U})\|_{L^2(\Omega_t^+)},
\end{equation}
\begin{equation}\label{el1'}
\quad \|\dot{V}(t)\|_{L^2(\mathbb{R}^3_-)}\leq \|(\dot{V},\p_t\dot{V})\|_{L^2(\Omega_t^-)},
\end{equation}
\begin{equation}\label{el2'}
\|\varphi (t)\|_{L^2(\Gamma )}\leq \|(\varphi,\p_t\varphi )\|_{L^2(\Gamma_t)}, \end{equation}
and
\begin{equation}\label{el2}
\|\nabla'\varphi (t)\|_{L^2(\Gamma )}\leq \|(\nabla'\varphi,\p_t\nabla'\varphi )\|_{L^2(\Gamma_t)} .
\end{equation}
Moreover, it follows from \eqref{qvntr} and the first boundary condition in \eqref{ELP1c} that
\begin{equation}\label{phit'}
\|\p_t\varphi (t)\|^2_{L^2(\Gamma )}\lesssim_K \|(\varphi,\nabla '\varphi )(t)\|^2_{L^2(\Gamma )} +\|f(t)\|^2_{L^2(\mathbb{R}^3_+)}+\nt \dot{U}(t)\nt^2_{H^1_{\rm tan}(\mathbb{R}^3_+)}.
\end{equation}
By applying \eqref{el1} to the source term $f$ we deduce from \eqref{phit'} that
\begin{equation}\label{phit}
\|\p_t\varphi (t)\|^2_{L^2(\Gamma )}\lesssim_K \|(\varphi,\nabla '\varphi )(t)\|^2_{L^2(\Gamma )} +\|f\|^2_{H^{1}_{\rm tan}(\Omega_T^+)}+\nt \dot{U}(t)\nt^2_{H^1_{\rm tan}(\mathbb{R}^3_+)}.
\end{equation}
In view of estimate \eqref{phit'} (integrated over the time interval), inequality \eqref{el2'}  yields
\begin{equation}\label{phi}
\|\varphi (t)\|_{L^2(\Gamma )}\lesssim_K \|f\|^2_{L^2(\Omega_T^+)}+\|(\varphi,\nabla '\varphi )\|^2_{L^2(\Gamma_t )}
+\|\dot{U}\|^2_{H^1_{\rm tan}(\Omega_t^+)}.
\end{equation}

\subsection{Closing the estimate in $H^1$}

Differentiating problem \eqref{ELP1} with respect to  $t$, $x_2$, and $x_3$ and applying the arguments analogous to those in \eqref{107}--\eqref{QQQQ}, we get the following counterparts of the energy inequality \eqref{enegineq} for $(\p_{\ell}\mathring{U},\p_{\ell}\mathring{V},\p_{\ell}\varphi )$:
\begin{multline}
I_{\ell}(t)+\sum_{\ell = 0,2,3}\int_{\Gamma}\bigg\{\frac{\mathfrak{s}|\nabla'\p_{\ell}\varphi |^2}{|\mathring{N}|^3} + \mathring{\mu}\p_{\ell}\varphi \p_{\ell}\dot{E}_N\bigg\}  + \int_{\Gamma_t}\mathcal{L}_{\ell} \\
\lesssim_K \|f\|^2_{H^1_{\rm tan}(\Omega_T^+)} +\|\dot{U}\|^2_{H^1_{\rm tan}(\Omega_t^+)}+\|\dot{V}\|^2_{H^1(\Omega_t^-)}+\|(\varphi ,\p_{\ell}\varphi,\nabla'\p_{\ell}\varphi )\|^2_{L^2(\Gamma_t )}  ,\quad \ell =0,2,3,
\label{enegineq1}
\end{multline}
where $\p_0:=\p_t$,
\[
I_{\ell}(t)= \int_{\mathbb{R}^3_+}A_0(\mathring{U})\p_{\ell}\dot{U}\cdot \p_{\ell}\dot{U} +\int_{\mathbb{R}^3_-}\mathfrak{B}_0(\mathring{v})\p_{\ell}\dot{V}\cdot\p_{\ell}
\dot{V},
\]
and $\mathcal{L}_{\ell}$ is the sum of terms like ${\rm coeff}\,\p_{\ell}\dot{q}\p_{\ell}\varphi$,
${\rm coeff}\,\p_{\ell}\dot{v}_N\p_{\ell}\varphi$, ${\rm coeff}\,\p_{\ell}\dot{h}_i\p_{\ell}\varphi$, ${\rm coeff}\,\p_{\ell}\dot{E}_i\p_{\ell}\varphi$ (with $i=1,2,3$), ${\rm coeff}(\p_{\ell}\varphi)^2$,  ${\rm coeff}\,\varphi\p_{\ell}\dot{q}$,
${\rm coeff}\,\varphi\p_{\ell}\dot{v}_N$, etc.
Here, we in particular used the fact that $\p_{\ell}\widetilde{A}_1(\mathring{U},\mathring{\Phi})|_{\Gamma}=0$ (see \cite{CDS14,MT14}) implying the estimate
\[
\int_{\Omega_t^+} \big(\p_{\ell}\widetilde{A}_1(\mathring{U},\mathring{\Phi})\big)\p_1\dot{U}\cdot
\p_{\ell}\dot{U} \lesssim_K \|\sigma\p_1\dot{U}\|^2_{L^2(\Omega_t^+)} \lesssim_K
\|\dot{U}\|^2_{H^1_{\rm tan}(\Omega_t^+)}
\]
for the lower-order term $\big(\p_{\ell}\widetilde{A}_1(\mathring{U},\mathring{\Phi})\big)\p_1\dot{U}\cdot
\p_{\ell}\dot{U}$ appearing in the energy identity for $\p_{\ell}\dot{U}$ (from this energy identity we deduce the energy inequality \eqref{enegineq1}).

Let us first consider $\ell =2$ or $\ell =3$ in \eqref{enegineq1}. Integrating by parts and using the Young inequality, we obtain
\begin{multline}\label{del1}
\int_{\Gamma}\mathring{\mu}\p_{\ell}\varphi \p_{\ell}\dot{E}_N =
-\int_{\Gamma}\big\{\mathring{\mu} \dot{E}_N\p_{\ell}^2\varphi + (\p_{\ell}\mathring{\mu})\p_{\ell}\varphi \dot{E}_N\big\} \\ \lesssim_K \delta \|(\p_{\ell}\varphi ,\p_{\ell}^2\varphi )(t)\|^2_{L^2(\Gamma )} +
\frac{1}{\delta}\|\dot{E}_N(t)\|^2_{L^2(\Gamma )}\qquad
\end{multline}
for all $\delta >0$. Then, the crucial point is that instead of the usage of a standard trace theorem we pass to the volume integral and apply the Young inequality with $\delta^2$ for estimating the trace $\dot{E}_N|_{\Gamma}$:
\begin{equation}\label{del2}
\|\dot{E}_N(t)\|^2_{L^2(\Gamma )}=2\int_{\mathbb{R}^3_-}\dot{E}_N\p_1\dot{E}_N\leq
\delta^2\|\p_1\dot{E}_N(t)\|^2_{L^2(\mathbb{R}^3_-)}+\frac{1}{\delta^2}
\|\dot{E}_N(t)\|^2_{L^2(\mathbb{R}^3_-)}.
\end{equation}
By virtue of \eqref{d1V} and \eqref{el1'}, combining \eqref{del1} and \eqref{del2} gives
\begin{multline}\label{del3}
\int_{\Gamma}\mathring{\mu}\p_{\ell}\varphi \p_{\ell}\dot{E}_N \\ \lesssim_K
\delta \left(\|(\p_{\ell}\varphi ,\p_{\ell}^2\varphi )(t)\|^2_{L^2(\Gamma )} +
 \nt\dot{V}(t)\nt^2_{H^1(\mathbb{R}^3_-)}\right) +
\frac{1}{\delta^3}\|\dot{V}\|^2_{H^1(\Omega_t^-)},\quad \ell =2,3.
\end{multline}
Regarding the terms contained in the sums $\mathcal{L}_2$ and $\mathcal{L}_3$, they are estimated by integrating by parts and applying  inequalities \eqref{qvntr} and \eqref{Vtr}, for example,
\begin{multline*}
\int_{\Gamma_t}{\rm coeff}\,\p_2\dot{q}\p_2\varphi = -\int_{\Gamma_t}\big\{{\rm coeff}\,\dot{q}\p_2^2\varphi + \p_2({\rm coeff})\,\dot{q}\p_2\varphi\big\} \\
\lesssim_K
\|f\|^2_{H^1_{\rm tan}(\Omega_T^+)} +\|\dot{U}\|^2_{H^1_{\rm tan}(\Omega_t^+)}+\|\dot{V}\|^2_{H^1(\Omega_t^-)}+\|(\p_2\varphi,\p_2^2\varphi )\|^2_{L^2(\Gamma_t )}.
\end{multline*}
We finally get
\begin{equation}\label{L23}
\int_{\Gamma_t}\mathcal{L}_{\ell} \lesssim_K
\|f\|^2_{H^1_{\rm tan}(\Omega_T^+)} +\|\dot{U}\|^2_{H^1_{\rm tan}(\Omega_t^+)}+\|\dot{V}\|^2_{H^1(\Omega_t^-)}+\|(\varphi,\p_{\ell}\varphi,
\p_{\ell}^2\varphi )\|^2_{L^2(\Gamma_t )},\quad \ell =2,3.
\end{equation}

We now consider the case $\ell =0$  in \eqref{enegineq1}. In view of \eqref{div-}, from \eqref{ELP1b} we can deduce the following equation (see \cite{CDS14})
\begin{equation}\label{tEN}
\varepsilon\p_t\dot{E}_N -\p_2\big( \dot{h}_3+\dot{h}_1\p_2\mathring{\Phi}+\varepsilon\dot{E}_3\p_t\mathring{\Phi}\big)
+\p_3\big( \dot{h}_2+\dot{h}_1\p_3\mathring{\Phi}-\varepsilon\dot{E}_2\p_t\mathring{\Phi}\big)=0
\qquad  \textrm{in } \Omega^-_T.
\end{equation}
Expressing $\p_t\dot{E}_N$ from \eqref{tEN} through $x_2$- and $x_3$-derivatives, we reduce the integral
$
\int_{\Gamma}\mathring{\mu}\p_t\varphi \p_t\dot{E}_N
$
contained in \eqref{enegineq1} for $\ell =0$ to integrals like
\begin{equation}\label{tEN'}
\int_{\Gamma} {\rm coeff}\,\p_t\varphi\p_2\dot{h}_1\quad \mbox{or}\quad
\int_{\Gamma} {\rm coeff}\,\dot{h}_1\p_t\varphi .
\end{equation}
We estimate the first integral in \eqref{tEN'} by using \eqref{qvn} and the arguments analogous to those in \eqref{del1} and \eqref{del2}:
\begin{align}
\nonumber
\int_{\Gamma} {\rm coeff}\,\p_t\varphi\p_2\dot{h}_1 &= -\int_{\Gamma}\big\{ {\rm coeff}\,\dot{h}_1\p_t\p_2\varphi + \p_2({\rm coeff})\,\dot{h}_1\p_t\varphi \big\} \\
\nonumber
&
\lesssim_K \delta \|(\p_t\varphi ,\p_t\p_2\varphi )(t)\|^2_{L^2(\Gamma )} +
\frac{1}{\delta}\|\dot{h}_1(t)\|^2_{L^2(\Gamma )} \\
&
\lesssim_K \delta \left(\|(\p_t\varphi ,\p_t\p_2\varphi )(t)\|^2_{L^2(\Gamma )} +
 \nt \dot{U}(t)\nt^2_{H^1_{\rm tan}(\mathbb{R}^3_+)}\right) +
\frac{1}{\delta^3}\|\dot{U}(t)\|_{L^2(\mathbb{R}^3_-)}.
\label{del4}
\end{align}
Here and everywhere else the positive constant $\delta$ can be taken the same as in \eqref{del1}. Clearly, the second integral in \eqref{tEN'} is also estimated by the right-hand side in the second inequality in \eqref{del4}. We finally obtain the estimate
\[
\int_{\Gamma}\mathring{\mu}\p_t\varphi \p_t\dot{E}_N
\lesssim_K \delta \left(\|(\p_t\varphi ,\p_t\p_2\varphi )(t)\|^2_{L^2(\Gamma )} +
 \nt \dot{U}(t)\nt^2_{H^1_{\rm tan}(\mathbb{R}^3_+)}\right) +
\frac{1}{\delta^3}\|\dot{U}(t)\|_{L^2(\mathbb{R}^3_-)},
\]
which, in view of \eqref{el1} and \eqref{phit}, implies
\begin{multline}\label{del5}
\int_{\Gamma}\mathring{\mu}\p_t\varphi \p_t\dot{E}_N \\
\lesssim_K \delta \left(\|(\varphi ,\nabla'\varphi,\p_t\p_2\varphi )(t)\|^2_{L^2(\Gamma )} +
 \nt \dot{U}(t)\nt^2_{H^1_{\rm tan}(\mathbb{R}^3_+)}\right) +  \|f\|^2_{H^{1}_{\rm tan}(\Omega_T^+)}+
\frac{1}{\delta^3}\|\dot{U}\|^2_{H^1_{\rm tan}(\Omega_t^+)}.
\end{multline}

Regarding the terms contained in the sum $\mathcal{L}_0$, we first describe the estimation of the most ``dangerous'' terms  ${\rm coeff}\,\p_{t}\dot{q}\p_{t}\varphi$,
${\rm coeff}\,\p_{t}\dot{v}_N\p_{t}\varphi$, ${\rm coeff}\,\p_{t}\dot{h}_i\p_{t}\varphi$, ${\rm coeff}\,\p_{t}\dot{E}_i\p_{t}\varphi$ ($i=1,2,3$) whereas the rest terms are estimated in a simpler way (see a short comment below). Let us, for example, consider the term ${\rm coeff}\,\p_{t}\dot{h}_2\p_{t}\varphi$. The crucial role in its estimation is played by expressing $\p_{t}\varphi$ from the first boundary condition in \eqref{ELP1c} through $\dot{v}_N$, $\nabla'\varphi$, and $\varphi$. This reduces $\int_{\Gamma_t}{\rm coeff}\,\p_{t}\dot{h}_2\p_{t}\varphi$ to the integrals
\begin{equation}\label{int3}
\int_{\Gamma_t} {\rm coeff}\,\dot{v}_N\p_{t}\dot{h}_2 ,\quad \int_{\Gamma_t} {\rm coeff}\,\p_{t}\dot{h}_2\p_k\varphi \quad (k=2,3),\quad\mbox{and}\quad \int_{\Gamma_t} {\rm coeff}\,\varphi \p_{t}\dot{h}_2.
\end{equation}

We estimate the first integral in \eqref{int3} by integrating by parts and using the Young inequality (with $\delta$) together with \eqref{qvn}, \eqref{d1V}, and \eqref{el1}:
\begin{align*}
\int_{\Gamma_t} {\rm coeff} & \,\dot{v}_N\p_{t}\dot{h}_2 \\ & = \int_{\Omega_t^-} \big\{{\rm coeff}\,\p_1\dot{v}_N^-\p_{t}\dot{h}_2 +\p_1({\rm coeff})\dot{v}_N\p_{t}\dot{h}_2-\p_t({\rm coeff})\dot{v}_N^-\p_{1}\dot{h}_2
-{\rm coeff}\,\p_t\dot{v}_N^-\p_{1}\dot{h}_2\big\}\\ & \qquad -\int_{\mathbb{R}^3_-}
{\rm coeff}\,\dot{v}_N^-\p_{1}\dot{h}_2 \\ & \lesssim_K
\|f\|^2_{H^1_{\rm tan}(\Omega_T^+)} +\|\dot{U}\|^2_{H^1_{\rm tan}(\Omega_t^+)}+\|\dot{V}\|^2_{H^1(\Omega_t^-)} +\delta \|\p_1\dot{h}_2(t)\|^2_{L^2(\mathbb{R}^3_-)} +\frac{1}{\delta}\|\dot{v}_N(t)\|^2_{L^2(\mathbb{R}^3_+)} \\ &
\lesssim_K
\|f\|^2_{H^1_{\rm tan}(\Omega_T^+)} +\frac{1}{\delta}\|\dot{U}\|^2_{H^1_{\rm tan}(\Omega_t^+)}+\|\dot{V}\|^2_{H^1(\Omega_t^-)}+
\delta \nt\dot{V}(t)\nt^2_{H^1(\mathbb{R}^3_-)},
\end{align*}
where $\dot{v}_N^-(t,x):=\dot{v}_N^-(t,-x_1,x')$. The integral of the second type in \eqref{int3} is estimated as follows (let, for example, $k=2$):
\begin{align*}
\int_{\Gamma_t} {\rm coeff}\,\p_{t}\dot{h}_2\p_2\varphi  & = \int_\Gamma {\rm coeff}\,\dot{h}_2\p_2\varphi - \int_{\Gamma_t}\big\{ \p_{t}({\rm coeff})\dot{h}_2\p_2\varphi + {\rm coeff}\,\dot{h}_2\p_t\p_2\varphi\big\} \\
& \lesssim_K \frac{1}{\delta}\|(\nabla'\varphi,\p_t\nabla'\varphi )\|^2_{L^2(\Gamma_t)}  +\delta \|\dot{V} (t)\|^2_{L^2(\Gamma )} \\ & \qquad\quad +\|\dot{V}\|^2_{L^2(\Gamma_t )} +
\|(\p_2\varphi , \p_t\p_2\varphi )\|^2_{L^2(\Gamma_t )} \\ & \lesssim_K
\frac{1}{\delta}\|(\nabla'\varphi,\p_t\nabla'\varphi )\|^2_{L^2(\Gamma_t)}
+\delta \nt\dot{V}(t)\nt^2_{H^1(\mathbb{R}^3_-)} +\|\dot{V}\|^2_{H^1(\Omega_t^-)}
\end{align*}
Here we have used estimates \eqref{Vtr} and \eqref{el2}.

Clearly, the third integral in \eqref{int3} can be estimated in the same way as the second one. Moreover, the integrals of some terms contained in the sum $\mathcal{L}_0$ look like the third integral in \eqref{int3}, and the rest terms in $\mathcal{L}_0$ are even estimated in an easier way. We finally deduce the estimate
\begin{multline}\label{L0}
\int_{\Gamma_t}\mathcal{L}_0 \lesssim_K
\|f\|^2_{H^1_{\rm tan}(\Omega_T^+)} +\frac{1}{\delta}\left( \|(\varphi ,\nabla'\varphi,\p_t\nabla'\varphi )\|^2_{L^2(\Gamma_t)}+\|\dot{U}\|^2_{H^1_{\rm tan}(\Omega_t^+)}\right)+\|\dot{V}\|^2_{H^1(\Omega_t^-)} \\ +
\delta \left(\nt \dot{U}(t)\nt^2_{H^1_{\rm tan}(\mathbb{R}^3_+)}+\nt\dot{V}(t)\nt^2_{H^1(\mathbb{R}^3_-)}
\right).
\end{multline}

Taking into account  \eqref{del3}, \eqref{L23}, \eqref{del5}, and \eqref{L0} and combining \eqref{wder}, \eqref{d1V}, \eqref{el1}, \eqref{el1'}, \eqref{el2}, and \eqref{phit}--\eqref{enegineq1}, we can choose $\delta $ small enough to obtain
\begin{multline}\label{preen}
\nt \dot{U}(t)\nt^2_{H^1_{\rm tan}(\mathbb{R}^3_+)}+\nt\dot{V}(t)\nt^2_{H^1(\mathbb{R}^3_-)}
+\nt (\varphi ,\nabla'\varphi )(t)\nt^2_{H^1(\Gamma)} \\ \lesssim_K
\|f\|^2_{H^{1}_{\rm tan}(\Omega_T^+)}+
\|\dot{U}\|^2_{H^{1}_{\rm tan}(\Omega_t^+)}+\|\dot{V}\|^2_{H^{1}(\Omega_t^-)}+\|(\varphi , \nabla'\varphi )\|^2_{H^{1}(\Gamma_t)},
\end{multline}
where
\[
\nt (\cdot )(t)\nt^2_{H^1(\Gamma)}:= \| (\cdot )(t)\|^2_{H^1(\Gamma)}+
\| \p_t(\cdot )(t)\|^2_{L^2(\Gamma)}.
\]
By virtue of Gr\"{o}nwall's inequality, \eqref{preen} implies the desired a priori estimate \eqref{54}.

\section{Concluding remarks}
\label{sec:5}

We have verified the stabilizing effect of surface tension on the evolution of the free interface separating a perfectly conducting inviscid fluid governed by the ideal compressible MHD equations from a vacuum where the electric and magnetic fields satisfy the Maxwell equations. Indeed, if we do not take the influence of surface tension into account, then, as was shown in \cite{MT14}, the planar interface can be violently unstable in the sense of Kelvin--Helmholz-type instability, which is associated with the ill-posedness of the corresponding constant coefficient linearized problem.

At the same time, the proof of the local-in-time well-posedness of the original nonlinear free interface problem with nonzero surface tension is still an open problem. Actually, with the help of the a priori estimate \eqref{54} deduced for the linearized problem one can prove (by a standard argument, see, e.g., \cite{ST14}) the uniqueness of a smooth solution to the nonlinear problem. Moreover, we think that, having in hand the existence of smooth solutions to the linear problem \eqref{ELP1}, nonlinear existence can be achieved by a modified Nash--Moser iteration scheme (as, for example, in \cite{TW22a,TW22c}).

That is, the lack of the proof of the existence of solutions to the linearized problem is now the main obstacle towards the proof of the local well-posedness of the nonlinear problem \eqref{NP1}. Note that the classical duality argument \cite{CP}  cannot be directly applied to our linear problem \eqref{ELP1} because for it we are not able to close the a priori estimate in $L^2$.  A possible way of overcoming this difficulty could be the strategy in \cite{TW22a,TW22c} based on inventing a suitable regularized linear problem. We do not know yet how to realize this strategy in our case and postpone this work to the future research.

\end{document}